\documentstyle[11pt,mssymb]{article}

\newcommand{\arrow}{\longrightarrow}
\renewcommand{\phi}{\varphi}
\newcommand{\Z}{{\Bbb Z}}
\newcommand{\C}{{\Bbb C}}
\newcommand{\R}{{\Bbb R}}

\newcommand{\1}{\sqrt{-1}\:}
\newcommand{\g}{{\goth g}}
\newcommand{\inangles}[1]{{\langle #1\rangle}}

\newcommand{\inbfpare}[1]{{%
  \mbox{\tt (}\hspace{-5pt}\mbox{\tt (} #1 %
  \mbox{\tt )}\hspace{-5pt}\mbox{\tt )}%
}}

\makeatletter

   \renewcommand{\section}{%
     \@startsection
        {section}{0}{.6em}%
        {2.5\baselineskip}%
        {1\baselineskip}%
        {
          \setcounter{footnote}{0}\setcounter{equation}{0}%
         \large\bf}}

   \renewcommand{\thesection}{\arabic{section}.}

\newcounter{lemma}[section]
\renewcommand{\thelemma}{{Lemma \thesection\arabic{lemma}}}
\newcommand{\lemma}{%
     \refstepcounter{lemma}
     {\bf \thelemma:\ }}

\newcounter{sublemma}[section]
\newcounter{corollary}[section]
\renewcommand{\thecorollary}{{Corollary \thesection\arabic{corollary}}}
\newcommand{\corollary}{%
     \refstepcounter{corollary}
     {\bf \thecorollary:\ }}

\newcounter{theorem}[section]
\renewcommand{\thetheorem}{{Theorem \thesection\arabic{theorem}}}
\newcommand{\theorem}{%
     \refstepcounter{theorem}
     {\bf \thetheorem:\ }}

\newcounter{proposition}[section]
\renewcommand{\theproposition}
       {{Proposition \thesection\arabic{proposition}}}
\newcommand{\proposition}{%
     \refstepcounter{proposition}
     {\bf \theproposition:\ }}

\newcounter{definition}[section]
\renewcommand{\thedefinition}
       {{Definition \thesection\arabic{definition}}}
\newcommand{\definition}{%
     \refstepcounter{definition}
     {\bf \thedefinition:\ }}

\newcounter{example}[section]
\newcommand{\eqref}[1]{(\ref{#1})}

\newcommand{\ps@verbit}{%
  \renewcommand{\@oddhead}{%
          \scriptsize
          {Holomorphic symplectic geometry \  {\rm II}}
          \hfil\tiny {M. Verbitsky, February 1994}}
  \renewcommand{\@evenhead}{\@oddhead}
  \renewcommand{\@oddfoot}{\hfil\thepage\hfil}
  \renewcommand{\@evenfoot}{\@oddfoot}}

\makeatother

\def\blacksquare{\hbox{\vrule width 4pt height 4pt depth 0pt}}

\begin{document}

\centerline{\bf Hyperk\"ahler embeddings}
\centerline{\bf and holomorphic symplectic geometry \  {\rm II}.}

\centerline{Mikhail Verbitsky,}
\centerline{verbit@math.harvard.edu}

\hfill

\hfill

{\bf 0. Introduction.}

This is a second part of the treatment of
complex analytic subvarieties of a holomorphically symplectic
K\"ahler manifold.
For the convenience of the reader, in the first two sections of this
paper we recall the  definitions and
results of the first part (\cite{_part_one_}).

\hfill

By Calabi-Yau theorem, the holomorphically symplectic K\"ahler manifolds
can be supplied with a  Ricci-flat Riemannian metric.
This implies that such manifolds
are hyperk\"ahler (\ref{_hyperkaehler_manifold_Definition_}).
Conversely, all hyperk\"ahler manifolds are
holomorphically symplectic (Proposition 2.1).

For a closed analytic subvariety $S$ of a holomorphically symplectic
$M$, one can restrict the holomorphic symplectic form of $M$
to the Zarisky tangent sheaf of $S$. If this restriction
is non-degenerate outside of singularities of $S$
and the same is true for the set of singular points of $S$,
this subvariety is called non-degenerately symplectic
(\ref{_non-dege_symple_Definition_}).
Of course, non-degenerately symplectic
subvarieties are of even complex dimension.

The hyperk\"ahler manifold is endowed with the canonical
$SU(2)$-action in this cohomology space. Fix an induced
complex structure on a hyperk\"ahler manifold $M$.
Let $N$ be a closed analytic subset of $M$. Then $N$ defines
a cycle $[N]$ in cohomology of $M$. Denote the Poincare
dual cocycle by $\inangles{N}$.
In \cite{_part_one_} we proved that if $\inangles{N}$ is
$SU(2)$-invariant then $N$ is non-degenerately
symplectic (Theorem 2.2 of \cite{_part_one_}).

Take a generic element $M_0$ in a given deformation class of
compact holomorphically symplectic K\"ahler manifolds. Then
all integer $(p,p)$-cycles on $M_0$ (i. e.,
all elements of $H^{2p}(M,\Z)\cap H^{p,p}(M)$)
are $SU(2)$-invariant (\ref{_generic_are_dense_Proposition_}).
According to Theorem 2.2 of \cite{_part_one_},  this immediately implies the
following statement: All closed analytic subvarieties of $N$
are non-degenerately symplectic. If such subvariety is smooth,
it is also a hyperk\"ahler manifold
(\ref{_symplectic_=>_hyperkaehler_Proposition_}).

Let $M$ be a hyperk\"ahler
manifold with three complex structures $I$, $J$ and $K$.
The closed subset $X\in M$ is called {\bf tri-analytic}
if $X$ is analytic with respect to $I$, $J$ and $K$.

In this paper we prove the following.
Let $N$ be a closed analytic subset of the compact
holomorphic symplectic manifold $M$. Assume that
$[N]$ is $SU(2)$-invariant. Then $N$ is trianalytic
(\ref {_G_M_invariant_implies_trianalytic_Theorem_}).
In particular, all closed analytic submanifolds of
the compact holomorphic symplectic manifold
of generic type are endowed with the hyperk\"ahler
structure which is preserved by the embedding
(\ref{_hyperkae_embeddings_Corollary_}).

For the compact complex torus this statement
implies a following corollary: if $N\subset M$
is a closed analytic subvariety of a generic compact
torus, then $N$ is a point. However, this statement
seems to be well known.


\section{Hyperk\"ahler manifolds.}


\definition \label{_hyperkaehler_manifold_Definition_} 
(\cite{_Beauville_},
\cite{_Besse:Einst_Manifo_}) A {\bf hyperk\"ahler manifold} is a
Riemannian manifold $M$ endowed with three complex structures $I$, $J$
and $K$, such that the following holds.

\hspace{5mm}   (i)  $M$ is K\"ahler with respect to these structures and

\hspace{5mm}   (ii) $I$, $J$ and $K$, considered as  endomorphisms
of a real tangent bundle, satisfy the relation
$I\circ J=-J\circ I = K$.

\hfill

This means that the hyperk\"ahler manifold has the natural action of
quaternions ${\Bbb H}$ in its real tangent bundle.
Therefore its complex dimension is even.


Let $\mbox{ad}I$, $\mbox{ad}J$ and $\mbox{ad}K$ be the operators on the
bundles of differential forms over a hyperk\"ahler manifold
$M$ which are defined as follows. Define $\mbox{ad}I$.
Let this operator act as a complex structure operator
$I$ on the bundle of differential 1-forms. We
extend it on $i$-forms for arbitrary $i$ using Leibnitz
formula: $\mbox{ad}I(\alpha\wedge\beta)=\mbox{ad}I(\alpha)\wedge\beta+
\alpha\wedge \mbox{ad}I(\beta)$. Since Leibnitz
formula is true for a commutator in a Lie algebras, one can immediately
obtain the following identities, which are implied by the same
identities in ${\Bbb H}$:

\[
   [\mbox{ad}I,\mbox{ad}J]=2\mbox{ad}K;\;
   [\mbox{ad}J,\mbox{ad}K]=2\mbox{ad}I;\;
\]

\[
   [\mbox{ad}K,\mbox{ad}I]=2\mbox{ad}J
\]

Therefore, the operators $\mbox{ad}I,\mbox{ad}J,\mbox{ad}K$
generate a Lie algebra $\goth{su}(2)$ acting on the
bundle of differential forms. We can integrate this
Lie algebra action to the action of a Lie group
$G_M=SU(2)$. In particular, operators $I$, $J$
and $K$, which act on differential forms by the formula
$I(\alpha\wedge\beta)=I(\alpha)\wedge I(\beta)$,
belong to this group.

{\bf Proposition 1.1:} There is an action of the Lie group $SU(2)$
and Lie algebra $\goth{su}(2)$ on the bundle of differential
forms over a hyperk\"ahler manifold. This action is
parallel, and therefore it commutes with Laplace operator.
$\blacksquare$

If $M$ is compact, this implies that there is
a canonical $SU(2)$-action on $H^i(M,\R)$ (see
\cite{_so5_on_cohomo_}).

\hfill

Let $M$ be a hyperk\"ahler manifold with a Riemannian form $<\cdot,\cdot>$.
Let the form $\omega_I := <I(\cdot),\cdot>$ be the usual K\"ahler
form  which is closed and parallel
(with respect to the connection). Analogously defined
forms $\omega_J$ and $\omega_K$ are
also closed and parallel.

The simple linear algebraic
consideration (\cite{_Besse:Einst_Manifo_}) shows that \hfill
$\omega_J+\sqrt{-1}\omega_K$ is of
type $(2,0)$ and, being closed, this form is also holomorphic.
It is called {\bf the canonical holomorphic symplectic form
of a manifold M}. Conversely, if there is a parallel
holomorphic symplectic form on a K\"ahler manifold $M$,
then this manifold has a hyperk\"ahler structure (\cite{_Besse:Einst_Manifo_}).

If some $compact$ K\"ahler manifold $M$ admits non-degenerate
holomorphic symplectic form $\Omega$, the Calabi-Yau
(\cite{_Yau:Calabi-Yau_}) theorem
implies that $M$ is hyperk\"ahler (Proposition 2.1).
This follows from the existence of a K\"ahler
metric on $M$ such that $\Omega$ is parallel for the Levi-Civitta
connection associated with this metric.

\hfill

Let $M$ be a hyperk\"ahler manifold with complex structures
$I$, $J$ and $K$. For any real numbers $a$, $b$, $c$
such that $a^2+b^2+c^2=1$ the operator $L:=aI+bJ+cK$ is also
an almost complex structure: $L^2=-1$.
Clearly, $L$ is parallel with respect to connection.
This implies that $L$ is a complex structure, and
that $M$ is K\"ahler with respect to $L$.

\hfill

\definition \label{_induced_structures_Definion_}
If $M$ is a  hyperk\"ahler manifold,
the complex structure $L$ is called {\bf induced
by a hyperk\"ahler structure}, if $L=aI+bJ+cK$ for some
real numbers $a,b,c\:|\:a^2+b^2+c^2=1$.

\hfill

\hfill

If $M$ is a hyperk\"ahler manifold and $L$ is induced complex structure,
we will denote $M$, considered as a complex manifold with respect to
$L$, by $(M,L)$ or, sometimes, by $M_L$.

\hfill

Consider the Lie algebra $\goth{g}_M$ generated by ${ad}L$ for all $L$
induced by a hyperk\"ahler structure on $M$. One can easily see
that $\goth{g}_M=\goth{su}(2)$.
The Lie algebra $\goth{g}_M$ is called {\bf isotropy algebra} of $M$, and
corresponding Lie group $G_M$ is called an {\bf isotropy group}
of $M$. By Proposition 1.1, the action of the group is parallel,
and therefore it commutes with Laplace operator in differential
forms. In particular, this implies that the action of the isotropy
group $G_M$ preserves harmonic forms, and therefore this
group canonically acts on cohomology of $M$.

\hfill

\proposition \label{_G_M_invariant_forms_Proposition_} 
Let $\omega$ be a differential form over
a hyperk\"ahler manifold $M$. The form $\omega$ is $G_M$-invariant
if and only if it is of Hodge type $(p,p)$ with respect to all
induced complex structures on $M$.

{\bf Proof:} Assume that $\omega$ is $G_M$-invariant.
This implies that all elements of $\g_M$ act trivially on
$\omega$ and, in particular, that $\mbox{ad}L(\omega)=0$
for any induced complex structure $L$. On the other hand,
$\mbox{ad}L(\omega)=(p-q)\1$ if $\omega$ is of Hodge type $(p,q)$.
Therefore $\omega$ is of Hodge type $(p,p)$ with respect to any
induced complex structure $L$.

Conversely, assume that $\omega$ is of type $(p,p)$ with respect
to all induced $L$. Then $\mbox{ad}L(\omega)=0$ for any induced $L$.
By definition, $\g_M$ is generated by such $\mbox{ad}L(\omega)=0$,
and therefore $\g_M$ and $G_M$ act trivially on $\omega$. $\blacksquare$


\section{Holomorphic symplectic geometry.}


\definition \label{_holomorphi_symple_Definition_} 
The compact K\"ahler manifold $M$ is called
holomorphically symplectic if there is a holomorphic 2-form $\Omega$
over $M$ such that $\Omega^n=\Omega\wedge\Omega\wedge...$ is
a nowhere degenerate section of a canonical class of $M$.
There, $2n=dim_\C(M)$.

Note that we assumed compactness of $M$.%
\footnote{If one wants to define a holomorphic symplectic
structure in a situation when $M$ is not compact,
one should require also the equation $\nabla'\Omega$ to held.
The operator $\nabla':\;\Lambda^{p,0}(M)\arrow\Lambda^{p+1,0}(M)$
is a holomorphic differential defined on differential $(p,0)$-forms
(\cite{_Griffiths_Harris_}).}
One observes that the holomorphically symplectic
manifold has a trivial canonical bundle.
A hyperk\"ahler manifold is holomorphically symplectic
(see Section 1). There is a converse proposition:

\proposition \label{_symplectic_=>_hyperkaehler_Proposition_}
(\cite{_Beauville_}, \cite{_Besse:Einst_Manifo_})
Let $M$ be a holomorphically
symplectic K\"ahler manifold with the holomorphic symplectic form
$\Omega$, a K\"ahler class
$[\omega]\in H^{1,1}(M)$ and a complex structure $I$.
There is a unique hyperk\"ahler structure $(I,J,K,(\cdot,\cdot))$
over $M$ such that the cohomology class of the symplectic form
$\omega_I=(\cdot,I\cdot)$ is equal to $[\omega]$ and the
canonical symplectic form $\omega_J+\1\omega_K$ is
equal to $\Omega$.

\ref{_symplectic_=>_hyperkaehler_Proposition_} immediately
follows from the Calabi-Yau theorem (\cite{_Yau:Calabi-Yau_}).
$\:\;\blacksquare$

\hfill

For each complex analytic variety $X$ and a point $x\in X$,
we denote the Zariski tangent space to $X$ in $x$ by $T_xX$.

\definition \label{_non-dege_symple_Definition_}
Let $M$ be a holomorphically symplectic
manifold and $S\subset M$ be its closed complex analytic subvariety.
Assume that $S$ is closed in $M$ and reduced. It is called
{\bf non-degenerately symplectic} if for each point $s\in S$
outside of the singularities of $S$ the restriction
of the holomorphic symplectic form $\Omega$ to $T_sM$
is nondegenerate on $T_s S\subset T_s M$, and the set
$Sing(S)$ of the singular points of $S$ is nondegenerately
symplectic. This definition refers to itself, but
since $dim\:Sing(S)<dim\:S$, it is consistent.

Of course, the complex dimension of a non-degenerately symplectic
variety is even.

\hfill

Let $M$ be a holomorphically symplectic K\"ahler manifold.
By \ref{_symplectic_=>_hyperkaehler_Proposition_},
$M$ has a unique hyperk\"ahler metric with
the same K\"ahler class and holomorphic symplectic form.
Therefore one can without ambiguity speak
about the action of $G_M$ on $H^*(M,\R)$ (see Proposition 1.1).
Of course, this action essentially depends on the choice
of K\"ahler class.

\hfill

\definition \label{_generic_manifolds_Definition_} 
Let $\omega\in H^{1,1}(M)$ be the K\"ahler of a K\"ahler matric
defined on a holomorphically symplectic
manifold $M$. We say that $\omega$ {\bf induces the $SU(2)$-action
of general type} when all elements of the group
\[ H^{pp}(M)\cap H^{2p}(M,\Z)\]
are $G_M$-invariant. The
action of $SU(2)\cong G_M$ is defined by
\ref{_symplectic_=>_hyperkaehler_Proposition_}.
The holomorphically symplectic manifold $M$ is
called {\bf of general type} if there
exists a K\"ahler class on $M$ which induces
an $SU(2)$-action of general type.

\hfill

As \ref{_subvarieties_of_generic_mfold_are_nondege_symple_Theorem_}
implies, the holomorhically symplectic manifold
of general type has no Weil divisors. Therefore these manifolds
have connected Picard group. In particular, such manifolds
are never algebraic.

\hfill

\proposition \label{_generic_are_dense_Proposition_} 
Let $M$ be a hyperk\"ahler manifold. Let $S$
be the set of induced complex structures over $M$. Let $S_0\subset S$
be the set of $R\in S$ such that the natural K\"ahler metric
on $(M,R)$ induces the $SU(2)$ action of general type.
Then $S_0$ is dense in $S$.

{\bf Proof:}
Let $A$ be the set of all $\alpha\in H^{2p}(M,\Z)$ such that $\alpha$
is not $G_M$-invariant. The set $A$ is countable. For each
$\alpha\in A$, let $S_\alpha$ be the set of all $R\in S$
such that $\alpha$ is of type $(p,p)$ with respect to $R$.
The set $S_0$ of all induced complex structures of general type
is equal to $\{S\backslash\bigcup_{\alpha\in A}S_\alpha\}$.
Now, to prove \ref{_generic_are_dense_Proposition_}
it is sufficient to show that
$S_\alpha$ is a finite set for each $\alpha\in A$.
This would imply that $S_0$ is a complement of a
countable set to a 2-sphere $S$, and therefore dense in $S$.

As it follows from Section 1,
$\alpha$ is of type $(p,p)$ with respect to $R$ if and
only if $ad \:R(\alpha)=0$. Now, let $V$ be a representation
of $\goth{su}(2)$, and $v\in V$ be a non-invariant vector.
It is easy to see that the element $a\in \goth{su}(2)$
such that $a(v)=0$ is unique up to a constant, if it exists.
This implies that if $\alpha$ is not $G_M$-invariant
there are no more than two $R\in S$ such that
$ad R(\alpha)=0$.
Of course, these two elements
of $S$ are opposite to each other. $\;\blacksquare$

\hfill

One can easily deduce from the results in
\cite{_Todorov:Moduli_I_II_} and
from \ref{_generic_are_dense_Proposition_} that the set of points
associated with holomorphically symplectic
manifolds of general type is dense in the classifying space
of holomorphically symplectic manifolds.

\hfill

For a K\"ahler manifold $M$ , $m=dim_\C M$
and a form $\alpha\in H^{2p}(M,\C)$, define
\[deg(\alpha):=\int_M L^{m-p}(\alpha)\]
where $L$ is a Hodge operator of exterior
multiplication by the K\"ahler form $\omega$
(see \cite{_Griffiths_Harris_}).
Of course, the degree of forms of Hodge type
$(p,q)$ with $p\neq q$ is equal zero,
so only $(p,p)$-form can possibly have
non-zero degree.

\hfill

We recall that the real dimension of a
holomorphically symplectic manifold is divisible by 4.

\theorem \label{_G_M_invariant_cycles_over_Theorem_} 
(Theorem 2.1 of \cite{_part_one_}).
Let $M$ be a holomorhically symplectic K\"ahler manifold with
a holomorphic symplectic form $\Omega$.
Let $\alpha$ be a $G_M$-invariant form of non-zero degree.
Then the dimension of $\alpha$ is divisible by 4. Moreover,
\[ \int_M \Omega^n\wedge\bar\Omega^n\wedge\alpha=2^n deg(\alpha),\]
where $n=\frac{1}{4}(dim_\R M-dim\:\alpha)$.

$\blacksquare$

\hfill

\proposition \label{_subvarieties_of_gene_are_even-dimen_} %
(Theorem 2.2 of \cite{_part_one_}). Let $M$ be a holomorphic symplectic
manifold of general type. All closed analytic subvarieties
of $M$ have even complex dimension.

{\bf Proof:} \ref{_subvarieties_of_gene_are_even-dimen_} immediately
follows from \ref{_G_M_invariant_cycles_over_Theorem_}.
$\;\;\blacksquare$

\hfill

\theorem \label{_subvarieties_of_generic_mfold_are_nondege_symple_Theorem_} %
(Theorem 2.3 of \cite{_part_one_}).
Let $M$ be a holomorphic symplectic
manifold of general type. All reduced closed analytic
subvarieties of $M$ are non-degenerately symplectic.

$\blacksquare$

\hfill

Combining this with \ref{_symplectic_=>_hyperkaehler_Proposition_},
one obtains

\theorem 
Let $M$ be a holomorphically symplectic
manifold of general type, and $S\subset M$ be its smooth differentiable
submanifold. If $S$ is analytic in $M$, it is a
hyperk\"ahler manifold. $\;\;\blacksquare$

\hfill

						       %
\section{When analytic implies tri-analytic}           %
						       %

\label{_analytic_implies_trianalytic_Section}

Let $M$ be a compact hyperk\"ahler manifold, $dim_\R M =2m$.
Let $I$ be an induced complex structure. As usually, $(M,I)$
denotes $M$ considered as a K\"ahler manifold with the
complex structure defined by $I$.

\hfill

\definition\label{_tri-analytic_Definition_} 
Let $N\subset M$ be a closed subset of $M$. Then $N$ is
called {\bf tri-analytic} if $N$ is an analytic subset
of $(M,I)$ for any induced complex structure $I$.

\hfill

Let $N\subset(M,I)$ be
a closed analytic subvariety of $(M,I)$, $dim_\C N= n$.
Let $[N]\in H_{2n}(M)$ denote the homology class
represented by $N$. Let $\inangles N\in H^{2m-2n}(M)$ denote
the Poincare dual cohomology class. Recall that
the hyperk\"ahler structure induces the action of
the group $G_M=SU(2)$ on the space $H^{2m-2n}(M)$.
The main result of this section is following:

\hfill

\theorem\label{_G_M_invariant_implies_trianalytic_Theorem_} 
Assume that $\inangles N\in  H^{2m-2n}(M)$ is invariant with respect
to the action of $G_M$ on $H^{2m-2n}(M)$. Then $N$ is tri-analytic%
\footnote{The number $n=dim_\C N$ is even by
\ref{_G_M_invariant_cycles_over_Theorem_}.}%
.%

\hfill

\ref{_G_M_invariant_implies_trianalytic_Theorem_} has the following
important corollary:

\corollary \label{_hyperkae_embeddings_Corollary_} 
Let $M$ be a holomorphically symplectic
manifold of general type, and $S\subset M$ be its smooth complex
submanifold. Let $\omega$ be the K\"ahler class
which induces an $SU(2)$-action of general type.
Pick a hyperk\"ahler metric $s$ associated with
$\omega$ by \ref{_symplectic_=>_hyperkaehler_Proposition_}.
Then the restriction of $s$ to $S$ is a hyperk\"ahler
metric.

$\;\;\blacksquare$

\hfill

The rest of this section is dedicated to the proof
of \ref{_G_M_invariant_implies_trianalytic_Theorem_}.

\hfill

{\bf Proof of \ref{_G_M_invariant_implies_trianalytic_Theorem_}:}
We use the following proposition of linear algebra.
Let $V$ be a complex space and $V_\R$ be its underlying $\R$-linear
space. Let $W\subset V_\R$ be an $\R$-linear subspace of $V_\R$,
$dim\; W=2n$. Denote the Hermitian form on $V$ by
$\inbfpare {\cdot,\cdot}$. The $\R$-linear space $V_\R$ is
equipped with the positively defined symmetric scalar product
$(u,v):=Re\:\inbfpare {u,v}$. and the non-degenerated
symplectic form $\inangles{u,v}:= I\!m\,\inbfpare {u,v}$.
Let $(\cdot,\cdot)_W$ and $\omega=\inangles{\cdot,\cdot}_W$
be the restrictions of these forms to $W\subset V_\R$.
Clearly, $(\cdot,\cdot)_W$ is a positively defined
scalar product on $W$; therefore, $(\cdot,\cdot)_W$
in non-degenerate. By a scalar product on an $R$-vector
space, we define a {\bf volume form} as follows.

\hfill

\definition 
Let $H$ be an $\R$-linear space equipped with a
positively defined scalar product. Let $h=dim\; H$. The exterrior
form $V\!ol\in \Lambda^h(H)$ is called {\bf a volume form}
if the the standard hypercube with the side 1
has the volume 1 in the measure defined by $V\!ol$.

\hfill

The proof of correctness of this definition can be found
in any linear algebra textbook.
Clearly, the volume form is defined up to a sign.
This sign is determined by the choice of orientation on $H$.
In the same manner
we define the top degree differential form called
{\bf a volume form} $V\!ol$ on any oriented
Riemannian manifold.

\hfill

Let $V\!ol$ be the volume form on $W$ defined by the scalar product.
Then $V\!ol$ is a non-zero element of the 1-dimensional linear
space

\[
   \mbox{\bf Vol}=\Lambda^{2n}(W).
\]
The form $\omega^n=\inangles{\cdot,\cdot}_W^n$ is another
element of $\mbox{\bf Vol}$. The number $\frac{\omega^n}{Vol}\in \R$
is defined up to a sign, because $V\!ol$ is defined up to a sign.
Let $\eta_W:= |\frac{\omega^n}{Vol}|$. This number is an
invariant of $W$, $V$ and $\inbfpare{\cdot,\cdot}$.

\hfill

Consider the complex structure operator on $V$ as the
real endomorphism of $V_\R$:

\[ I:\;V_\R\arrow V_\R,\;\;\;
   I^2=-1.
\]

\hfill

\proposition\label{_Wirtinger_Proposition_} 
(Wirtinger's inequality)
Let $2n=dim_\R W$. Then $\eta_W\leq 2^n.$
Moreover, if $\eta_W=2^n$, then $I(W)=W$. In
other words, if $\eta_W=2^n$, then $W$ is a
complex subspace of $V$.

{\bf Proof:} \cite{_Stolzenberg_} page 7. $\;\;\blacksquare$

\hfill

We return to the situanion when $M$ is a hyperk\"ahler manifold
and $N\subset (M,I)$ is a closed analytic subvariety of $(M,I)$.
Let $J$ be an induced complex structure, and
$N^0 \subset M$ be the set of non-singular
points of $N$. Denote the standard embedding
$N^0 \hookrightarrow M$ by $\phi$.
Let $\omega_J\in \Lambda^{1,1}(M)$ be the K\"ahler
form induced by $J$. Denote the standard coupling of homology
and cohomology by

\[ \inangles{\cdot,\cdot}:\;
   H_i(M)\times H^i(M)\arrow \C.
\]

\hfill

\proposition\label{_restriction_is_coupling_Proposition_} 

\[ \inangles{[N],\omega_J^n}=
   \int\limits_{N^0}\phi^*(\omega_J^n).
\]

{\bf Proof:} Clear. $\;\;\blacksquare$

\hfill

Let $\inbfpare{\cdot,\cdot}_J$ denote the Hermitian
form on $M$ induced by the Riemannian metric and the
complex structure $J$. Let $\omega_J=\inangles{\cdot,\cdot}_J$
and $(\cdot,\cdot)$ denote its imaginary and real parts respectively.
It is clear that the scalar product $(\cdot,\cdot)$ is the original
Riemannian form on $M$. Thus, the real part of
$\inbfpare{\cdot,\cdot}_J$ does not depend on
the complex structure $J$.

Let $(\cdot,\cdot)_N:=\phi^*((\cdot,\cdot))$ be the
Riemannian form on $N^0$. Since $N$ is oriented,
$(\cdot,\cdot)_N$ defines a volume form $V\!ol\in \Lambda^{2n}(N^0)$.
By definition, $V\!ol$ is a nowhere degenerate section of the
1-dimensional\footnote {over $\R$}
 vector bundle $\Lambda^{2n}(N^0)$.
Since $N$ is analytic in $(M,I)$, we have
$V\!ol=1/2^n\phi^*(\omega_I)^n\in\Lambda^{2n}(N^0).$
Let $x\in N^0$. Consider $V=T_xM$ as a complex
space with the complex structure induced by $J$.
Let $W=T_xN^0$. Then $W\subset V_\R$ defines a number
$\eta_W$ as in \ref{_Wirtinger_Proposition_}. This
number depends on $x$ and $J$. For any induced complex
structure $J$, we define a function
$\eta_J:\; N^0\arrow \R^{\geq 0}$ which supplies the number
$\eta_W\in\R^{\geq 0}$ by the point $x\in N^0$.

\hfill

\proposition\label{_N_is_analytic_if_eta_is_constant_Proposition_} 
Let $J$ be an induced complex structure.
The closed set $N\subset M$ is analytic with respect to $J$
if and only if

\[ \forall x\in N^0 \;\;\;\; \eta_J(x)=2^n. \]

{\bf Proof:} The implication

\hfill

\centerline{($N$ is analytic w. r. to $J$)
            $\;\;\;\Rightarrow\;\;\;$ ($\eta_J(x)\equiv 2^n$)}

\hfill

\hspace{-1.5em}%
is clear because if $N$ is analytic with respect to $J$, then
$J(T_xN)=T_xN$ and $\eta_J(x)\equiv 2^n$ by
\ref{_Wirtinger_Proposition_}. We proceed proving
the converse implication. Assume that
$\forall x\in N^0$ we have $\eta_J(x)=2^n$. Then
$J(T_xN)=T_xN$ by \ref{_Wirtinger_Proposition_}.
Using Newlander-Nierenberg theorem, we see that
$N^0$ is an analytic subset of $(M,J)$.
Clearly, $N$ is a closure of $N^0$.
Since the closure of an analytic set is
also analytic, the set $N\subset M$ is also an
analytic subset of $(M,J)$. $\;\;\blacksquare$

\hfill

By definition of $\eta_J(x)$, we have

\begin{equation}\label{_eta_in_terms_of_omega_Equation_}
   \int\limits_{N^0}\phi^*(\omega_J^{2n})=
   \int\limits_{N^0} V\!ol\cdot \eta_J(x).
\end{equation}
By \ref{_Wirtinger_Proposition_}, we have $\eta_J(x)\leq 2^n$.
Clearly, the function
\[ \eta_J(x):\; N^0 \arrow \R^{\geq 0}\]
is continous. Therefore

\[ \int\limits_{N^0} V\!ol\cdot \eta_J(x)
   = 2^n \int\limits_{N^0} V\!ol
\]
if and only if $\eta_J(x)=2^n$ for every $x\in N^0$.
Combining this with
\eqref{_eta_in_terms_of_omega_Equation_}
and
\ref{_N_is_analytic_if_eta_is_constant_Proposition_},
we obtain the following statement:

\hfill

\proposition\label{_analiticity_in_terms_of_integrals_Proposition_} 

\nopagebreak

\vspace{\baselineskip}

\hspace{15mm}$\displaystyle\int\limits_{N^0}\phi^*(\omega_J^{2n})=
              2^n \int\limits_{N^0} V\!ol$

\hfill

\hspace{-1.9em}
if and only if $N$ is analytic with respect to $J$.

$\blacksquare$

\hfill

\hfill

Since $N$ is analytic with respect to $J$, we have

\[ 2^n \int\limits_{N^0} V\!ol = \int\limits_{N^0}\phi^*(\omega_J^{2n}). \]
Combining \ref{_analiticity_in_terms_of_integrals_Proposition_} with
\ref{_restriction_is_coupling_Proposition_}, we see that
\ref{_G_M_invariant_implies_trianalytic_Theorem_} is implied
by the following statement.

\hfill

\proposition\label{_volumes_of_G_M_invariant_cycles_Proposition_} 
Let $M$ be a compact hyperk\"ahler manifold, $dim_\C M=m$.
Let $[N]\in H_{2n}(M)$ be a homology class of $M$ such that
its Poincare dual cocycle $\inangles{N}\in H^{2m-2n}$ is
$G_M$-invariant. Let $I$ and $J$ be two induced complex
structures on $M$. Then $\inangles{[N], \omega_I^{2n}}=
\inangles{[N], \omega_J^{2n}}$.

{\bf Proof:} By definition,

\[ \inangles{[N], \omega_J^{2n}} =
   \int_M\inangles{N} \wedge \omega_J^{2n}.
\]
Therefore all we need to show is that
the number
\[ deg_J\alpha:= \int_M\alpha\wedge  \omega_J^2n\]
is independent on the choice of $J$ once the cohomology
class $\alpha$ is $G_M$-invariant.

Let $L_J:\; \Lambda^i(M)\arrow \Lambda^{i+2}(M)$
denote the Hodge operator acting on differential forms
over $M$, $L_J(\eta)=\omega_J\wedge \eta$.
Let $\Lambda_J:\; \Lambda^i(M)\arrow \Lambda^{i-2}(M)$
denote the adjoint operator. It is well known that
$L_J$, $\Lambda_J$ map harmonic form to harmonic ones.
Therefore one can consider $L_J$, $\Lambda_J$
as operators on the cohomology space $H^*(M)$.
Let $\goth a_M\subset End(H^*(M))$ be the Lie algebra
generated by $L_J$, $\Lambda_J$ for all induced complex
structures $J$. Let $\g_M\cong \goth{so}(3)$ be
the Lie algebra of $G_M$. Since $G_M$ non-trivially
acts on $H^*(M)$, we may consider $\g_M$ as a subalgebra
of $End(H^*(M))$. It is known that
$\goth a_M\cong \goth{so}(5)$ and that
$\g_M$ considered as a Lie subalgebra of $End(H^*(M))$
lies in $\goth a_M\subset End(H^*(M))$
(see \cite{_so5_on_cohomo_}).

\hfill

Let $H^*(M)=\oplus_{l\in\Pi}H_l$ be the isotypic
decomposition of a $\goth a_M$-module
$H^*(M)$. We recall that {\bf isotypic decomposition}
of a representation of an arbitrary reductive Lie algebra
is defined as follows. For each $l\in \Pi$, where
$\Pi$ is a weight lattice of $\goth{a}_M$,
the module $H_l$ is a union of all simple $\goth{a}_M$-submodules
of $H^*(M)$ with a highest weight $l$.
One can easily see that the isotypic decomposition does not depend
on a choice of a Cartan subalgebra of $\goth{a}_M$. This follows,
for example, from Schuhr's lemma.

\hfill

Let $H_o$ be the $\goth a_M$-submodule of $H^*(M)$ generated by
$H^0(M)\cong \C$. Clearly, $H_o$ is an isotypic component of
$H^*(M)$ (see \cite{_part_one_} for details). Let $\alpha_o$ be the component
of
$\alpha$ which corresponds to the summand
$H_o\subset \oplus_{l\in\Pi}H_l$. Since $\g_M\subset \goth a_M$,
the isotypic component $\alpha_o$ of $\alpha$ is $\g_M$-invariant.
It was proven that $deg_J(\alpha_o)=deg_J(\alpha)$ for all
induced complex structures $J$ (see the paragraph right
after the proof of Lemma 3.1 in \cite{_part_one_}).

Let $I$, $J$, $K$ be a triple of induced complex structures
on $M$, such that

\[ I\circ J=-J\circ I= K. \]
Lemma 3.2 of \cite{_part_one_} implies that

\[ \alpha_o=c(L_I^2+ L_J^2+L_K^2)^{\frac{m-n}{2}} {\Bbb I} \]
where $ {\Bbb I} $ is a generator of $H^0(M)$ and
$c$ is a constant. In this notation, the equality

\[ deg_I(\alpha_o)= deg_J(\alpha_o)= deg_K(\alpha_o) \]
is obvious. We proved the following lemma

\hfill

\lemma \label{_deg_is_independent_on_I_J_K_Lemma_} 
If $I$, $J$, $K$ are induced complex structures
on $M$, such that

\[ I\circ J=-J\circ I= K, \]
then

\[ deg_I(\alpha_o)= deg_J(\alpha_o)= deg_K(\alpha_o). \]

$\blacksquare$

\hfill

We deduce \ref{_volumes_of_G_M_invariant_cycles_Proposition_}
from \ref{_deg_is_independent_on_I_J_K_Lemma_} as follows.
Let $V$ be the space of purely imaginary quaternions.
By definition, the space of purely imaginary quaternions
$V\subset \Bbb H$ is a three-dimensional vector space over $\R$':

\[ V= \{ t\in {\Bbb H} \; | \; \bar t =-t \}, \]

\[ V= \{ aI+bJ+cK,\; a,b,c \in \R\}. \]

The set of induced complex structures can be considered
as a sphere $S$ of radius 1 in $V$. For $I$, $J$ in $S$,
we have $I\circ J=-J\circ I$ if and only if the vector
$\stackrel\arrow {(0,I)}\in V$ is perpendicular to
$\stackrel\arrow {(0,J)}\in V$. On the other hand,
if $\stackrel\arrow {(0,I)}\bot\stackrel\arrow {(0,J)}$,
then $K:= I\circ J$ belongs to $S$ and the vectors
$\stackrel\arrow {(0,I)}$, $\stackrel\arrow {(0,J)}$,
$\stackrel\arrow {(0,K)}$ are pairwise orthogonal.
We obtain that the triples
$I$, $J$, $K$ of induced complex structures which satisfy

\[ I\circ J=-J\circ I= K \]
are in one-to-one correspondence with the orthonormal repers in
$V$. Therefore \ref{_deg_is_independent_on_I_J_K_Lemma_}
implies the following statement:

\hfill

\lemma \label{_orthogonal_quaternions_and_degree_Lemma_} 
Let $I$, $J$ be the induced complex structures such
that the vectors $\stackrel\arrow {(0,I)}\in V$ and
$\stackrel\arrow {(0,J)}\in V$ are orthogonal. Then
$deg_I(\alpha)=deg_J(\alpha)$.

$\blacksquare$

\hfill

\ref{_orthogonal_quaternions_and_degree_Lemma_} trivially
implies \ref{_volumes_of_G_M_invariant_cycles_Proposition_}.
\ref{_G_M_invariant_implies_trianalytic_Theorem_} is proven.


\hfill

{\bf Acknowledgements:}
I am very grateful to my advisor David Kazhdan for a warm support
and encouragement. Professor Y.-T. Siu was extremely helpful
kindly answering my questions. He provided me with the
reference to the Stolzenberg's book (\cite{_Stolzenberg_}).
F. A. Bogomolov cleared some of the misconceptions about
the holomorphic symplectic geometry that I had.
I am grateful to Roman Bezrukavnikov for interesting
discussions and insightful remarks.
I am also grateful to MIT math department
for allowing me the use of their computing
facilities.

\hfill

\end{document}